\documentclass{gtmon_a}
\pdfoutput=1

\usepackage{amscd}
\usepackage{amsxtra}
\usepackage[cmtip,arrow]{xy}   
\usepackage{pb-diagram,pb-xy}  

\proceedingstitle{Proceedings of the Nishida Fest (Kinosaki 2003)}
\conferencestart{28 July 2003}
\conferenceend{8 August 2003}
\conferencename{International Conference in Homotopy Theory}
\conferencelocation{Kinosaki, Japan}

\editor{Matthew Ando}
\givenname{Matthew}
\surname{Ando}

\editor{Norihiko Minami}
\givenname{Norihiko}
\surname{Minami}

\editor{Jack Morava}
\givenname{Jack}
\surname{Morava}

\editor{W Stephen Wilson}
\givenname{W Stephen}
\surname{Wilson}

\title{Nilpotency of the Bauer--Furuta stable homotopy\\Seiberg--Witten invariants}

\author{Mikio Furuta}
\givenname{Mikio}
\surname{Furuta}
\address{Department of Mathematical Sciences\\
University of Tokyo\\\newline
Tokyo 153-8914\\
Japan}
\email{furuta@ms.u-tokyo.ac.jp}
\urladdr{}

\author{Yukio Kametani}
\givenname{Yukio}
\surname{Kametani}
\address{Department of Mathematics\\
Keio University\\\newline
Yokohama 223-8522\\
Japan}
\email{kametani@math.keio.ac.jp}
\urladdr{}

\author{Norihiko Minami}
\givenname{Norihiko}
\surname{Minami}
\address{Nagoya Institute of Technology\\\newline
Nagoya 466-8555\\
Japan}
\email{minami.norihiko@nitech.ac.jp}
\urladdr{}

\dedicatory{Dedicated to Professor Goro Nishida for his 60th birthday}

\volumenumber{10}
\issuenumber{}
\publicationyear{2007}
\papernumber{7}
\startpage{147}
\endpage{154}

\doi{}
\MR{}
\Zbl{}

\arxivreference{}


\received{1 January 2005}
\revised{5 January 2006}
\accepted{}
\published{29 January 2007}
\publishedonline{29 January 2007}
\proposed{}
\seconded{}
\corresponding{}
\version{}

\makeatletter
\def\cnewtheorem#1[#2]#3{\newtheorem{#1}{#3}[section]
\expandafter\let\csname c@#1\endcsname\c@theorem}

\let\xysavmatrix\xymatrix
\def\xymatrix{\disablesubscriptcorrection\xysavmatrix}
\AtBeginDocument{\let\bar\wbar\let\tilde\wtilde\let\hat\what}

\newcommand{\Spin}{\mathrm{Spin}}

\newtheorem{theorem}{Theorem}[section]
\cnewtheorem{lemma}[theorem]{Lemma}
\cnewtheorem{proposition}[theorem]{Proposition}
\cnewtheorem{corollary}[theorem]{Corollary}
\cnewtheorem{question}[theorem]{Question}

\theoremstyle{definition}
\cnewtheorem{definition}[theorem]{Definition}
\cnewtheorem{example}[theorem]{Example}
\cnewtheorem{xca}[theorem]{Exercise}

\theoremstyle{remark}
\cnewtheorem{answer}[theorem]{Answer}
\cnewtheorem{hint}[theorem]{Hint}
\cnewtheorem{remark}[theorem]{Remark}
\cnewtheorem{notation}[theorem]{Notation}

\numberwithin{equation}{section}

\makeatother  

\newcommand{\M}{{ \mathcal{M} }}

\newcommand{\CP}{\mathbb{CP}}

\makeop{proper}
\makeop{Map}
\makeop{id}
\makeop{sign}
\makeop{Ker}
\makeop{Res}
\makeop{Image}

\makeop{level}
\makeop{slevel}

\begin{document}

\begin{asciiabstract} 
We prove a nilpotency theorem for the Bauer-Furuta 
stable homotopy
Seiberg-Witten invariants for smooth closed
4-manifolds with trivial first Betti number.
\end{asciiabstract}

\begin{abstract} 
We prove a nilpotency theorem for the Bauer--Furuta 
stable homotopy
Seiberg--Witten invariants for smooth closed
4--manifolds with trivial first Betti number.
\end{abstract}

\maketitle

\section{Introduction}

As a byproduct of his celebrated work \cite{F1} on
the Yukio Matsumoto $11/8$ birthday conjecture \cite{Mat},
the first author and S Bauer \cite{Bau2,BF} significantly
refined the Seiberg--Witten invariants \cite{Wit} for
smooth closed $\Spin^c\ 4$--manifolds.
While the connected sum of two smooth connected $\Spin^c\ 4$--manifolds 
with $b_2^+ \geq 2$ always yields trivial
Seiberg--Witten invariants, these refined Seiberg--Witten
invariants yielded non-trivial values for some connected sums.

In this paper, we investigate the behavior of these refined
Seiberg--Witten invariants, which we call the Bauer--Furuta
stable homotopy Seiberg--Witten invariants,
with respect to the connected sum operations.
For simplicity, we always assume $b_1=0$.
Then our main theorem is the following:

\begin{theorem} \label{MT}
Let $X$ 
be a $\Spin^c\ 4$--manifold with
$b_2^+(X) \geq 1$ and $b_1(X)=0.$  
Then, there is some large $N$ such
that, for any $n\geq N$, $\#^n X$, the $n$--fold connected 
sum of $X$ 
with itself, has a
trivial Bauer--Furuta stable homotopy Seiberg--Witten 
invariant for any $\Spin^c$--structure $c$ of $\#^n X$
and any orientation $o_{\#^n X}$ of $H^+(\#^n X)$.
\end{theorem}

This is an easy consequence of
the Devinatz--Hopkins--Smith nilpotency theorem
\cite{DHS,HS}, which generalizes the earlier work 
of Nishida's nilpotency theorem \cite{Nis}.
Actually, the gluing theorem of S Bauer \cite{Bau1}
translates the connected sum operation corresponding
to the $S^1$--equivariant join operation, in which
case the usual non-equivariant smash product may 
be used.

After recalling the basic notations and backgrounds
in \fullref{sec2}, we present the Bauer gluing theorem 
in \fullref{sec3} and
the Devinatz--Hopkins--Smith nilpotency theorem
in \fullref{sec4}. Then, in \fullref{sec5},
we prove the main theorem.


We would like to express our sincere gratitudes to
the editor Matthew Ando and the referee, for their
useful suggestions and careful readings of the original
manuscript.  The third author was partially supported by Grant-in-Aid
for Scientific Research No.13440020, Japan Society for the
Promotion of Science.

The nilpotency is really the key philosophy of this paper, and we would
like to dedicate this paper to Professor Goro Nishida.

\section{Notation and background} \label{sec2}

We first set the following notations:
\begin{center}
\begin{tabular}{ll}
$X$        & an oriented closed 4--manifold with $b_1(X)=0$ \\
$c$        & a $\Spin^c$--structure of $X$\\
$\sign(X)$ & the signature of $X$ \\
$H^+(X)$   & the maximal positive definite subspace of $H^2(X,\R)$ \\
$o_X$      & an orientation of $H^+(X)$ \\
$b^+_2(X)$ & the dimension of $H^+(X)$
\end{tabular}
\end{center}
We now put
$$m := \tfrac18\bigl(c_1(L_c)^2- \sign(X)\bigr) \qquad
n :=b^+_2(X)$$
and assume $m\geq 0$. 
Actually,
we can always make $m 
\geq 0$ by changing the orientation of $X$,
if necessary.
Now the Bauer--Furuta stable homotopy Seiberg--Witten invariant
$\widetilde{SW}(X,c,o_X)$ \cite{BF} is defined for the data 
$(X,c,o_X)$ so that
\begin{enumerate}
\renewcommand{\labelenumi}{(\roman{enumi})}
\item \label{sec2.i} $SW(X,c,o_X) \in \{S(\C^{m}), S({\R}^n)\}^{U(1)} :=$

\hfill$\lim_{p,q\rightarrow \infty} 
          [S(\C^{p+m} \oplus {\R}^{q}),
           S(\C^{p}\oplus {\R}^{q+n})]^{U(1)}$
\item \label{sec2.ii} (via the $G$--Freudenthal
  suspension theorem by Hauschild \cite{Hau}, 
  Namboodiri \cite{Nam}) 
  If $n=b^+_2(X) \geq 2,$
  \begin{multline*}
     \{S(\C^{m}), S({\R}^n)\}^{U(1)} =
        \lim_{p,q\rightarrow \infty}
        [S(\C^{p+m} \oplus {\R}^{q}),
        S(\C^{p}\oplus {\R}^{q+n})]^{U(1)}  \\
       \overset{\cong}{\leftarrow}
        \lim_{q\rightarrow \infty}
        [S(\C^{m} \oplus {\R}^{q}),
        S({\R}^{q+n})]^{U(1)}  \\
      =  [ \CP^{m-1}*S^{q-1}, S^{n-1}*S^{q-1} ]
      =  \pi^{n-1}_s( \CP^{m-1}_+ ),
  \end{multline*}
  where $\dim\M = 2m-n-1$.
\end{enumerate}
Under the cohomotopy Hurewicz map %
\[
  \pi^{n-1}_s( \CP^{m-1}_+ ) 
  \to H^{n-1}( \CP^{m-1}_+ ),
\]
$SW(X,c,o_X)$  maps to the usual Seiberg--Witten invariant. 
%
%
%
%
%
%
Note that the difference of the top cell dimensions
of  $\CP^{m-1}_+$ and $S^{n-1}$ is nothing but the 
dimension of the
Seiberg--Witten moduli space $\M$:
\[
  \dim \M =  2( m - 1 ) - ( n -1 ) = 
2 m - n - 1.
\]

\section{Bauer's gluing theorem} \label{sec3}

For oriented closed 4--manifolds $X_1, X_2,$ 
consider their connected sum, the oriented closed 4--manifold
$X_1\# X_2.$ Then any $\Spin^c$--structure $c$ of
$X_1\# X_2$ and any orientation $o_{X_1\# X_2}$ of 
$H^+(X_1\# X_2),$ the maximal positive definite subspace 
of $H^2(X_1\# X_2,\R),$ are both joins of the corresponding
quantities of $X_1$ and $X_2.$ More explicitly, for
$i=1,2,$ there are $\Spin^c$--structures $c_i$ of
$X_i$ and orientations $o_{X_i}$ of 
$H^+(X_i),$ the maximal positive definite subspace 
of $H^2(X_i,\R),$ such that
$$c = c_1 \# c_2  \qquad
o_{X_1\# X_2} = o_{X_1}\# o_{X_2}$$
Thus, we may write
$$(X_1\# X_2,c,o_{X_1\# X_2}) =
(X_1\# X_2,c_1 \# c_2,o_{X_1}\# o_{X_2})
= (X_1,c_1,o_{X_1})\# (X_2,c_2,o_{X_2}).$$
Now the Bauer gluing formula \cite{Bau1}
describes the Bauer--Furuta stable homotopy Seiberg--Witten
invariant of this connected sum:

\begin{theorem}
The connected sum $\#$ corresponds to the join $*$:
$$SW( (X_1,c_1,o_{X_1})\#(X_2,c_2,o_{X_2}))
= SW(X_1,c_1,o_{X_1})*SW(X_2,c_2,o_{X_2}),$$
where the join
$$*\co \{S(\C^{2{k_1}}), S({\R}^{l_1})\}^{U(1)} \times
   \{S(\C^{2{k_1}}), S({\R}^{l_1})\}^{U(1)}  \\
\to \{S(\C^{2{k_1+k_2}}), S({\R}^{l_1+l_2})\}^{U(1)}$$
is induced by the usual join operation 
$S(U)*S(V) = S(U\oplus V).$
\end{theorem}

We note the dimension of the moduli space
behaves as
\begin{multline*}
\dim\M( (X_1,c_1,o_{X_1})\#(X_2,c_2,o_{X_2}))  \\
= \dim\M(X_1,c_1,o_{X_1})+\dim\M(X_2,c_2,o_{X_2})+1,
\end{multline*}
via the join operation.  From this, we 
observe the usual Seiberg--Witten invariant of
\[
(X_1,c_1,o_{X_1})\#(X_2,c_2,o_{X_2})
\]
is $0,$ when $b^+(X_i) \geq 2,\ (i=1,2),$
However, as was first pointed out by Bauer \cite{Bau1},
this is \emph{not} the case for the Bauer--Furuta stable homotopy
Seiberg--Witten invariant.  Actually, most applications
of the Bauer--Furuta stable homotopy Seiberg--Witten
invariants reflect this phenomenon.

\section{Nilpotency} \label{sec4}

We now change gears to hard-core homotopy theory.
We recall the fundamental concept and the
theorem of Hopkins--Smith \cite{HS}.

\begin{definition}{\rm \cite[Definition 1.i]{HS}}\qua
A map of spectra
\( 
f\co F \to X
\)
is {\it smash nilpotent} if for  $n \gg 0$ the 
induced map $n$--fold smash products
\(
f^{(n)}\co F^{(n)} \to X^{(n)} 
\)
is null.
\end{definition}

\begin{theorem}{\rm \cite[Theorem 3.iii]{HS}}\qua \label{HS}
A map $f\co F \to X$ from a finite spectrum to a $p$--local
spectrum is smash nilpotent if and only if $K(n)_*f = 0$
for all $0\leq n\leq \infty.$
\end{theorem}

This deep theorem generalizes their earlier work
\cite{DHS} with Devinatz, and specializes to Nishida's
nilpotency on the stable homotopy groups of spheres
\cite{Nis} when $F$ and $X$ are both 
spheres of distinct dimensions.

\section{Proof} \label{sec5}

Before the proof, we first consider some general
properties of the Bauer--Furuta stable homotopy 
Seiberg--Witten invariant 
$SW( \#_{i=1}^n X_i, c, o_{\#_{i=1}^n X_i})$
of $\#_{i=1}^n X_i.$ 

First, as was discussed at the beginning of \fullref{sec3}, 
we may always write
\[
SW( \#_{i=1}^n X_i, c, o_{\#_{i=1}^n X_i})
= SW( \#_{i=1}^n (X_i,c_i,o_{X_i}) ) 
= *_{i=1}^n SW( X_i,c_i,o_{X_i} )
\]
for some choices of $\Spin^c$--structure $c_i$ and 
orientation $o_{X_i}$ of $H^+(X_i),$
the maximal positive definite subspace of 
$H^2(X_i,\R).$

For each $i = 1,\ldots, n,$ let $f_i$ be a 
representative of 
\[
SW( X_i,c_i,o_{X_i} )
\in \{S(\C^{a_i}), S({\R}^{b_i})\}^{U(1)} 
\]
where
$$a_i = \tfrac18\bigl(c_1(L_{c_i})^2 - \sign(X_i)\bigr),\qquad
b_i = b_2^{+}(X_i).$$
To simplify our notations,
we express every map ``unstably'' in this
section.  For instance, we write
\[
f_i\co S(\C^{a_i}) \to S({\R}^{b_i}),
\]
though the general precise expression is of the form
$f_i\co S(\C^{a_i}\oplus \R^q) \to S({\R}^{b_i+q})$
(cf \fullref{sec2}\eqref{sec2.i}).
However, such simplified notations shall not cause
any problem here.

When $b_2^+(X_i) \geq 2$ for each $i = 1,\ldots, n,$ 
$f_i$ admits a non-equivariant counterpart
as in \fullref{sec2}\ref{sec2.ii}.  We denote it also
``unstably'' as
\[
\tilde{f}_i\co \CP^{a_i - 1} \to S^{b_i - 1},
\]
though the general precise expression is of the form
\[
\tilde{f}_i\co \CP^{a_i - 1}*S^{q-1} 
\to S^{b_i - 1}*S^{q-1} 
\eqno{\text{(cf \fullref{sec2}\ref{sec2.ii})}}
\]
Then a non-equivariant counterpart $\wtilde f$
of a representative 
$\bigwedge_{i=1}^n \wtilde{f}_i$ of
\[
SW( \#_{i=1}^n X_i, c, o_{\#_{i=1}^n X_i})
= SW( \#_{i=1}^n (X_i,c_i,o_{X_i}) ) 
= *_{i=1}^n SW( X_i,c_i,o_{X_i} )
\]
shows up at the top row of the following
commutative diagram (which is written
``unstably'' to simplify the notations):
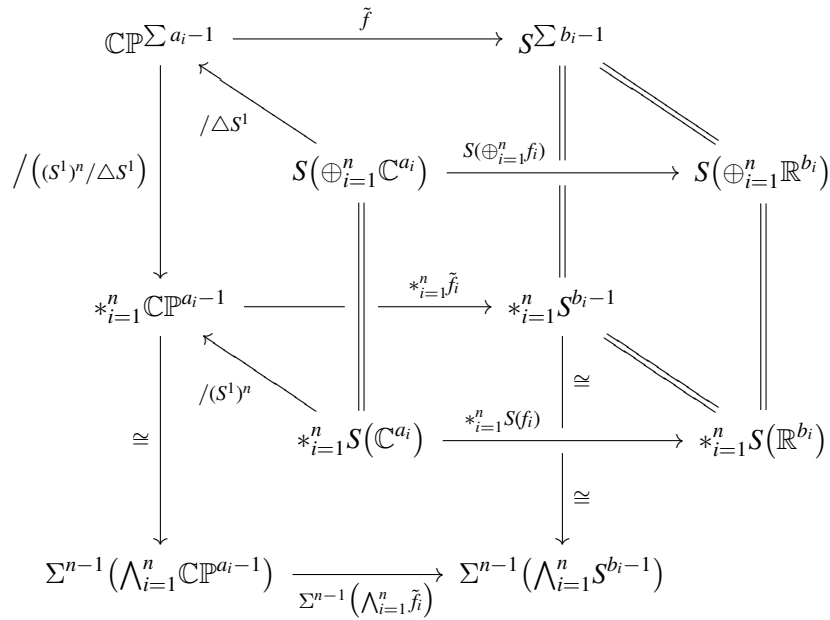
\[
\begin{diagram}
 \node{\CP^{\sum a_i - 1}} 
   \arrow[2]{e,t}{\wtilde f}
   \arrow[2]{s,l}{\big/\bigl( (S^1)^n/\triangle S^1 \bigr)}
 \node[2]{S^{\sum b_i - 1}}
   \arrow{se,=}{}
   \arrow{s,=}{}    
                       \\
 \node[2]{S\bigl(\oplus_{i=1}^n \C^{a_i}\bigr) }
   \arrow[2]{e,t,1}{S(\oplus_{i=1}^n f_i) }
   \arrow{nw,r}{/\triangle S^1}
   \arrow[2]{s,=}{}
 \node{} \arrow{s,=}{}
 \node{S\bigl(\oplus_{i=1}^n \R^{b_i}\bigr)}
   \arrow[2]{s,=}{} 
                       \\
 \node{*_{i=1}^n \CP^{ a_i - 1}}
   \arrow{e,-}{}
   \arrow[2]{s,l}{\cong}
 \node{}
   \arrow{e,t}{*_{i=1}^n \wtilde{f}_i }
 \node{*_{i=1}^n S^{b_i - 1}}
   \arrow{se,=}{} 
   \arrow{s,r,-}{\cong}
                       \\
 \node[2]{*_{i=1}^n S\bigl(\C^{a_i}\bigr) }
   \arrow[2]{e,t,1}{*_{i=1}^n S( f_i) }
   \arrow{nw,r}{/(S^1)^n}
 \node{}
   \arrow{s,r}{\cong}
 \node{*_{i=1}^n S\bigl(\R^{b_i}\bigr)}
                         \\
 \node{\Sigma^{n-1}
 \bigl({\textstyle\bigwedge_{i=1}^n} \CP^{ a_i - 1} \bigr)}
   \arrow[2]{e,b}{\Sigma^{n-1}
   \bigl( \bigwedge_{i=1}^n \wtilde{f}_i \bigr)}
 \node[2]{\Sigma^{n-1}
   \bigl( {\textstyle\bigwedge_{i=1}^n} S^{ b_i - 1} \bigr)}
\end{diagram}
\]


From this commutative diagram, we 
see immediately that, to show the
Bauer--Furuta stable homotopy Seiberg--Witten
invariant
\[
SW( \#_{i=1}^N X_i, c, o_{\#_{i=1}^N X_i})
= SW( \#_{i=1}^N (X_i,c_i,o_{X_i}) ) 
= *_{i=1}^N SW( X_i,c_i,o_{X_i} )
\]
is trivial, it is enough to show that the bottom row
\[
\Sigma^{n-1}\bigl( {\textstyle\bigwedge_{i=1}^n} \tilde{f}_i \bigr)\co
\Sigma^{n-1}
 \bigl( {\textstyle\bigwedge_{i=1}^n} \CP^{ a_i - 1} \bigr) \to
\Sigma^{n-1}
   \bigl( {\textstyle\bigwedge_{i=1}^n} S^{ b_i - 1} \bigr)
\]
is non-equivariantly (and, of course, stably)
trivial. 

Now we are ready to prove the theorem.

\begin{proof}[Proof of \fullref{MT}]
Since $b_2^{+}(X) \geq 1,$ we see that
$b_2^{+}(X\# X) \geq 2.$
Thus we may assume $b_2^{+}(X) \geq 2$ from the
beginning. Then we may apply the preceding
consideration with $X_i = X$ for $i=1,\ldots,n.$

Now the crucial, but very elementary, observation
is that there are just finitely many choices of
a $\Spin^c$--structure $c$ of $X$ and 
an orientation $o_{X}$ of $H^+(X).$ So, list all
such structures as $(X,c_j,o_{X_j})\ (j=1,\ldots,s).$
Now, for any $m,$ choose $n$ so that
$n\geq (m-1)s+1.$  Then, up to permutations
of factors, the Bauer--Furuta stable
homotopy Seiberg--Witten invariant
\[
SW( \#_{i=1}^n X_i, c, o_{\#_{i=1}^n X_i})
\]
of an arbitrary $\Spin^c$--structure $c$ of
$\#_{i=1}^n X_i$ and arbitrary orientation
$o_{\#_{i=1}^n X_i}$ of $H^+(o_{\#_{i=1}^n X_i})$
may be expressed in the form
\[
(*^m f_j)*g
\]
for some equivariant map $g$ and a representative
$f_j$ of $SW(X, c_j,o_{X_j})$ for some
$j\in \{1,\ldots,s\}.$

Thus, it remains to show any 
$S^1$--equivariant map
\[
f_0\co S(\C^a) \to S(\R^b)
\]
with $a\geq 1, b\geq 2$
is $S^1$--equivariantly join nilpotent (join
stably). Actually, if this is shown,
there are some $n_j$ so that $*^{n_j} f_j$
is $S^1$--equivariantly (join stably) nilpotent,
and so we may simply take 
\[
N = (\max\{n_j,j=1,\ldots,s\} - 1)s +1
\]
as $N$ in the statement of \fullref{MT}.

To show $f_0\co S(\C^a) \to S(\R^b)$ with $a\geq 1, b\geq 2$
is  $S^1$--equivariantly join nilpotent (join stably), 
we may assume its non-equivariant counter part 
$\tilde f_0\co \CP^{a-1} \to S^{b-1}$ induces the
trivial maps of the Morava $K$--theories $K(n)_*$
for all $n\geq 0.$ This is because joining once will
reduce to the case $b$ even.

Thus, $\wtilde f_0\co \CP^{a-1} \to S^{b-1}$ is
(non-equivariantly) smash-nilpotent (stably)
thanks to the Hopkins--Smith \fullref{HS}.

Now by the commutative diagram at the beginning
of this section with $a_i = a,\ b_i = b$ for any $i,$
the claim follows. 
\end{proof}

\textbf{Note added in proof}\qua
To show the smash-nilpotency of 
$f_0,$ we can also use more familiar $MU_*$--theory,
rather than the less familiar Morava $K$--theories 
$K(n)_*\ (n\geq 0).$ 
Actually, for $h_*= MU_*\ \text{or}\ K(n)_*,$ applying
the Atiyah--Hirzebruch spectral sequence, we see the 
groups $h_*(X)$ are concentrated in even dimensions for X a 
complex projective space and concentrated in odd dimensions 
for an odd dimensional sphere. We then use the Hopkins--Smith 
\fullref{HS} in the $MU_*$ version.

\bibliographystyle{gtart}
\bibliography{link}

\end{document}